\documentclass[10pt]{amsart}

\usepackage{amssymb,url,mathrsfs,euscript}
\usepackage[british]{babel}

\theoremstyle{remark}     
\newtheorem{rmk}{Remark}[section]\newtheorem*{rmk*}{Remark}

\newtheorem*{ex*}{Example}

\newtheorem*{acknowledgements}{Acknowledgements}

\swapnumbers              
\theoremstyle{plain}      
\newtheorem{lemma}[rmk]{Lemma}
\newtheorem{proposition}[rmk]{Proposition}
\newtheorem{theorem}[rmk]{Theorem}
\newtheorem{corollary}[rmk]{Corollary}

\theoremstyle{definition} 
\newtheorem{definition}[rmk]{Definition}

\numberwithin{equation}{section}

\newcommand{\hodge}{{*}}
\newcommand{\lie}[1]{\mathfrak{#1}}
\newcommand{\Lie}[1]{\textsl{#1}}
\newcommand{\G}{G\sb 2}
\DeclareMathOperator{\SO}{\Lie{SO}}

\DeclareMathOperator{\SU}{\Lie{SU}}
\DeclareMathOperator{\U}{\Lie{U}}
\DeclareMathOperator{\su}{\lie{su}}

\newcommand{\g}{\lie{g}}
\newcommand{\n}{\lie{n}}
\newcommand{\s}{\lie{s}}
\newcommand{\h}{\lie{h}}

\DeclareMathOperator{\codim}{\textsl{codim}\,}
\newcommand{\W}{\mathcal{W}}
\newcommand{\X}{\mathcal{X}}

\newcommand{\psp}{\psi^+}
\newcommand{\psm}{\psi^-}
\newcommand{\pspm}{\psi^\pm}

\newcommand{\iso}{\cong}

\newcommand{\lto}{\longrightarrow}

\newcommand{\tto}{\mapsto}
\newcommand{\lan}{\langle}
\newcommand{\ran}{\rangle}

\newcommand{\rcomp}[1]{\mathopen{\big[\mkern-5mu\big[} #1
\mathclose{\big]\mkern-5mu\big]}}
\newcommand{\rreal}[1]{\bigl[#1\bigr]}

\newcommand{\p}[1]{\frac{\partial #1}{\partial \t}}

\newcommand{\oo}{\omega^2}
\newcommand{\1}{e^1}\newcommand{\2}{e^2}\newcommand{\3}{e^3}
\newcommand{\4}{e^4}\newcommand{\5}{e^5}\newcommand{\6}{e^6}
\newcommand{\7}{e^7}
\renewcommand{\geq}{\geqslant}\renewcommand{\leq}{\leqslant}

\newcommand{\ad}{\textsl{ad}}
\newcommand{\tr}{\textsl{tr}}

\newcommand{\q}{\quad}\newcommand{\qq}{\qquad}

\newcommand{\bproof}{\begin{proof}}
\newcommand{\eproof}{\end{proof}}

\newcommand{\bsub}{\begin{subequations}}
     \newcommand{\esub}{\end{subequations}}
\newcommand{\ba}{\begin{array}}\newcommand{\ea}{\end{array}}

\newcommand{\f}{\varphi}\newcommand{\ff}{\hodge\varphi}
\newcommand{\R}{\mathbb{R}}\newcommand{\C}{\mathbb{C}}

\newcommand{\om}{\omega}

\newcommand{\op}{\oplus}

\newcommand{\w}{\wedge}
\newcommand{\set}{\subseteq}

\newcommand{\hd}{\hat d}

\newcommand{\vs}{\vphantom{\int_W^M}}
\renewcommand{\t}{\textrm{\footnotesize T}}

\begin{document}

\title[Conformally parallel $\G$ structures on a class of
solvmanifolds]{Conformally parallel $\boldsymbol{\G}$ structures\\ on a
    class of solvmanifolds}
\author{Simon G.~Chiossi} \address[S.G.Chiossi]{Institut for Matematik og
     Datalogi, Syddansk Universitet, Campusvej 55, 5230 Odense M, Denmark}
\email{chiossi@imada.sdu.dk}
\author{Anna Fino} \address[A.Fino]{Dipartimento
     di Matematica, Universit\`a di Torino, via Carlo Alberto 10, 10123
     Torino, Italy} \email{fino@dm.unito.it}

\begin{abstract}
  Starting from a $6$-dimensional nilpotent Lie group $N$ endowed with an
  invariant $\SU(3)$ structure, we construct a homogeneous conformally
  parallel $G_2$-metric on an associated solvmanifold. We classify all
  half-flat $\SU(3)$ structures that endow the rank-one solvable extension
  of $N$ with a conformally parallel $G_2$ structure. By suitably 
  deforming the $\SU(3)$ structures obtained, we are able to describe the
  corresponding non-homogeneous Ricci-flat metrics with holonomy contained
  in $G_2$.  In the process we also find a new metric with exceptional
  holonomy.
\end{abstract}

\maketitle

\section{Introduction}

\noindent
A seven-dimensional Riemannian manifold $(Y, g)$ is called a $G_2$-manifold
if it admits a reduction of the structure group of the tangent bundle to
the exceptional Lie group $G_2$.  The presence of a $G_2$ structure is
equivalent to the existence of a certain type of three-form $\varphi$ on the
manifold. Whenever this $3$-form is covariantly constant with respect to
the Levi--Civita connection then the holonomy group is contained in $G_2$,
and the corresponding manifold is called parallel. The development of the
theory of explicit metrics with holonomy $G_2$ follows the by-now-classical
line of Bonan~\cite{Bonan:G2-Spin7}, Fern\'andez and
Gray~\cite{Fernandez-G:G2}, Bryant~\cite{Bryant:exceptional} and
Salamon~\cite{Bryant-S:exceptional}.  We shall review a few relevant facts
in section $2$.

Interesting non-compact examples are provided by Gibbons, L\" u, Pope,
Stelle in~\cite{Gibbons-LPS:domain-walls}, where incomplete Ricci-flat
metrics of holonomy $G_2$ with a $2$-step nilpotent isometry group $N$
acting on orbits of codimension one are presented. It turns out that these
metrics have scaling symmetries generated by a homothetic Killing vector
field, and are locally isometric (modulo a conformal change) to homogeneous
metrics on solvable Lie groups.  The solvable Lie group in question is
obtained by extending the isometry group of the original manifold, and can
be seen as the universal cover of the product of $\R$ with the $2$-step
nilmanifold corresponding to $N$, which is a compact quotient
$\Gamma\backslash N$ by a discrete uniform subgroup.  Solvmanifolds ---
that is solvable Lie groups endowed with a left-invariant metric --- and in
particular solvable extensions of nilpotent Lie groups provide instances of
homogeneous Einstein manifolds.  The fact that any nilpotent Lie algebra of
dimension $6$ admits a solvable extension carrying Einstein
metrics~\cite{Will:Einstein-7-solv} will be of the foremost importance.

We shall concentrate on {\it conformally parallel $G_2$ structures},
characterised by the fact that the Riemannian metric $g$ can be modified to
metric with holonomy a subgroup of $G_2$ by a transformation
\begin{displaymath}
      g\tto e^{2f}g,
\end{displaymath}
for some function $f$.

In the light of \cite{Will:Einstein-7-solv}, it is natural to study such
$G_2$ structures on a rank-one solvable extension of a metric
$6$-dimensional nilpotent Lie algebra $\n$ endowed with an $\SU(3)$
structure $(\omega, \psi^+)$ and a non-singular self-adjoint derivation $D$
which is diagonalisable by a unitary basis.  This last condition is
equivalent to $(DJ)^2 = (JD)^2$ and we show that this is the compatibility
that one has to impose between $D$ and the $\SU(3)$ structure in order to
obtain the non-compact examples found in~\cite{Gibbons-LPS:domain-walls}.

As shown in section \ref{sec:extension}, such an
extension is given by a metric Lie algebra $\s=\n \oplus \R H$ with bracket
\begin{displaymath}
     [H,U] = DU, \quad [U, V] = [U, V]_{\n\times\n},
\end{displaymath}
where $U, V\in\n$ and $H\perp\n, \|H\|=1$. The subscript denotes the Lie
bracket on $\mathfrak n$, and the inner product extends that of $\mathfrak
n$. There is a natural $G_2$ structure on the manifold $Y=N \times \R$
corresponding to the $3$-form
\begin{displaymath}
     \varphi=\omega\wedge H^\flat +\psp\in\Lambda^3T^*Y,
\end{displaymath}
where $\flat$ is the isomorphism of $T$ onto $T^*$ induced by the metric. The
Lie algebra $\s$ is isomorphic to each fibre of the principal fibration
$T^*Y\lto Y$, and we prove the

\medbreak \noindent
Main result. {\it $(Y, \varphi)$ is conformally parallel if and
    only if $\n$ is either $\R^6$, or $2$-step nilpotent but not isomorphic
    to the Lie algebra $\h_3 \oplus \h_3$,}

\medbreak\noindent where ${\mathfrak h}_3$ denotes the real $3$-dimensional
Heinsenberg algebra (cf.~\S \ref{sec:classification}).  The operator $D$
has the same eigenvalue type of the derivation considered by Will to
construct Einstein metrics on 7-dimensional solvmanifolds
\cite{Will:Einstein-7-solv}.

In section~\ref{sec:metrics} we describe explicitly the corresponding
metrics $g$ with holonomy a non-trivial subgroup of $G_2$. Half of such
metrics have ${\mathit {Hol}}(g) = G_2$ and stem from the three irreducible
$2$-step nilpotent Lie algebras. The remaining metrics have holonomy either
$\SU(2)$ or $\SU(3)$ and correspond to Lie algebras with abelian summands.
Using this we show that some metrics have also been considered by
\cite{Gibbons-LPS:domain-walls} in the study of special domain walls in
string theory.  We are able to produce a new metric
with holonomy equal to $G_2$, that arises from the $6$-dimensional Lie
algebra spanned by $\1,\ldots,\6$ with
\begin{displaymath}
    e_2=[e_5,e_4], \q e_3=[e_6,e_4]=[e_1,e_5]
\end{displaymath}
as the only non-trivial brackets.

The conformally parallel $G_2$ structure forces the initial $\SU(3)$
structure to be of a special kind, known in the literature as
half-flat~\cite{Chiossi-S:SUG}.  This turns out to be a useful notion,
which allows one to find explicit metrics with holonomy $G_2$ by
investigating the corresponding Hitchin flow~\cite{Hitchin:forms}.
Section~\ref{sec:evolution} is especially devoted to such a description. We
determine a solution of the evolutions equations and compare the resulting
$G_2$ holonomy metrics with the ones previously described. These rank-one
solvmanifolds $S$ admit then a pair of distinguished metrics.  The first is
the homogeneous Einstein metric with negative scalar curvature constructed
in \cite{Will:Einstein-7-solv}. The other arises by conformally changing a
homogeneous metric and possesses a homothetic Killing field, i.e.~a vector
field with respect to which the Lie derivative of $g$ is a multiple of the
identity; our investigation proves that it is also obtainable by evolving
the original $\SU(3)$ structure.

\begin{acknowledgements}
   The authors are indebted to S.~Salamon and A.~Swann for the invaluable
   suggestions, and thank I.~Agricola and T.~Friedrich for hospitality
   during the initial stage of this project.  They are both members of the
   \textsc{Edge} Research Training Network
   \textsc{hprn-ct-\oldstylenums{2000}-\oldstylenums{00101}}, supported by
   the European Human Potential Programme. The research is partially
   supported by \textsc{Miur, Gnsaga--Indam} in Italy.
\end{acknowledgements}

\section{$G$ structures in $6$ and $7$ dimensions}

\noindent
Suppose that $X$ indicates a six-dimensional nilmanifold with an invariant
almost Hermitian structure. Thus, $X$ is endowed with an orthogonal almost
complex structure $J$ and a non-degenerate $2$-form $\omega$ which induce a
Riemannian metric $h$. An $\SU(3)$-reduction of the structure group is
determined by fixing a real $3$-form $\psp$ lying in the $S^1$-bundle of
unit elements inside the canonical bundle $\rcomp{\Lambda^{3,0}}$ at each
point. We adopt the parenthetical notation of \cite{Salamon:holonomy} to
indicate real modules of $\{p,q\}$-forms underlying the complex space
$\Lambda^{p+q}_\C$. Let $\Psi=\psp+i\psm$ be the associated
holomorphic section  (so that $J\psm=-\psp$). The description is 
always intended to be local, so one
can define the forms
\begin{equation}
     \label{eq:standardstructure}
     \begin{gathered}
       \omega=e^{14}-e^{23}+e^{56},\\
       \psp+i\psm=(e^1+ie^4)\wedge(e^2-ie^3)\wedge(e^5+ie^6),
     \end{gathered}
\end{equation}
of type $(1,1)$ and $(3,0)$ relative to $J$. It has become customary to
suppress wedge signs when writing differential forms, so $e^{ij\dots}$
  indicates $e^i\wedge e^j\wedge\dots$ from now on.  Following
\cite{Chiossi-S:SUG}  and  \cite{Bor-HL:Bochner} we tackle
six-dimensional geometry by means of the enhanced Gray and
Hervella decomposition of the intrinsic torsion space into five
representations $\W_1,\ldots ,\W_5$. These are the $\SU(3)$-modules
appearing in $\Lambda^1\otimes \bigl(\rcomp{\Lambda^{2,0}}\oplus\R\bigr)$
that identify the kind of almost Hermitian structure. Complex
$\SU(3)$-manifolds are for instance characterised by the vanishing of the
intrinsic torsion components belonging to $\W_1\iso\R\oplus\R,
\W_2\iso\su(3)\oplus \su(3)$.  Tagging the irreducible `halves' by $\pm$,
one can correspondingly split the Nijenhius tensor $N_J=N^+_J+N^-_J$. The
modules $\W_1^\pm, \W_2^\pm$ can be defined explicitly by prescribing the
various types of the real forms
\begin{displaymath}
     \begin{array}{rcrcl}
       d\psp & = & -2W_5\wedge\psp+W_2^+\wedge\omega+W_1^+\oo \\
       d\psm & = & 2W_5\wedge\psm+W_2^-\wedge\omega+W_1^-\oo
     \end{array}
\end{displaymath}
corresponding to $\Lambda^4T^*X \iso \rcomp{\Lambda^{0,1}}\oplus
\rreal{\Lambda^{1,1}_0} \oplus \R$. The nought in the middle term denotes
$(1,1)$-forms $\alpha$ satisfying $\alpha\wedge\omega=0$, called primitive.

\smallbreak Moving up one dimension, we consider a product $Y$ of $X$ with
$\R$, endowed with metric $g$. Indicating by $\7$ the unit 1-form on the
real line one obtains a basis for the cotangent spaces $T_y^*Y$.  The
manifold $Y$ inherits a non-degenerate three-form $\f=\omega\wedge\7+\psp$
which is stable, {\sl \`a la} Hitchin~\cite{Hitchin:forms}, and defines a
reduction to the exceptional group.  The fundamental material for the $G_2$
story can be found in standard
references~\cite{Salamon:holonomy,Joyce:book}.  Let us only recall that the
Riemannian geometry of $Y$ is completely determined by the tensor
\begin{displaymath}
     \f=e^{125}-e^{345}+e^{567}+e^{136}+e^{246}-e^{237}+e^{147}.
\end{displaymath}
The seminal results of Fern\'andez and Gray~\cite{Fernandez-G:G2} permit
one to describe $G_2$ geometry exclusively in algebraic terms, by looking
at the various components of $d\f,d\ff$ in the irreducible summands
$\X_1,\X_2,\X_3$ and $\X_4$ of the space $T^*Y\otimes \g_2^\perp$. Many
authors have studied special classes of $G_2$ structures, see for
instance~\cite{Cabrera-MS:G2,Friedrich-I:skew,Cleyton-I:cG2}.  Before
concentrating on a particular situation, recall that in general the
exterior derivatives can be expressed as
\begin{displaymath}
     \left\{
       \begin{array}{rcl}
         d\ff & = & 4\tau_4\w\ff+\tau_2\w\f\\
         d\f  & = & \tau_1\w\ff+3\tau_4\w\f+\hodge\tau_3
       \end{array}
     \right.,
\end{displaymath}
where the various $\tau_i$'s represent the differential forms corresponding
to the representations $\X_i$, as in \cite{Bryant:G2}. For example $\tau_4$
is the $1$-form encoding the `conformal' data of the structure.  With the
convention of dropping all unnecessary wedge signs, the torsion three-form
of the unique $G_2$-connection~\cite{Friedrich-I:skew} is given
by
\begin{displaymath}
       \Phi=\tfrac 76\tau_1\f-\hodge d\f+\hodge(4\tau_4\,\f)
\end{displaymath}
in terms of $\tau_1=\tfrac 17 g( d\f,\ff )$ and $\tau_4=-\tfrac 34
\hodge(\hodge d\f\wedge\f)$, the latter being the Lee form of the
7-manifold, essentially.\\

Our aim is to study conformally parallel $G_2$ structures on Riemannian
products, otherwise said manifolds $X\times \R$ whose intrinsic torsion
belongs to the class $\X_4$ only. If this is the case, the above pair of
equations simplifies to
\begin{displaymath}
      d\ff = 4\tau_4\w\ff , \qq  d\f  = 3\tau_4\w\f
\end{displaymath}
and the obstruction to the reduction of the holonomy can be written as
$\Phi=\hodge(\tau_4\,\f)$, proportional to the Hodge dual of $d\f$.  Now
$\tau_4$ is a closed 1-form in the more general setting of
$G_2$T-structures, so as soon as one has $\dim H^1(Y,\R)=1$ (see
\eqref{(*)}), it will be natural to assume it is proportional to $\7$. So
let us rewrite those relations as
\begin{equation}
     \label{eq:m}
     \left\{
       \begin{array}{rcl}
         d\ff & = & 4m\7\wedge\ff\\
         d\f  & = & 3m\7\wedge\f
       \end{array}
     \right.,
\end{equation}
which also serve as a definition for the real constant $m$. To prevent the
holonomy of the metric $g$ from reducing to $G_2$, we implicitly assume
that $m$ does not vanish.\\
    
We shall next fit the geometric picture into the theory of Lie algebras,
and suppose $X$ is a nilpotent Lie group. This is indeed no real
restriction since~\cite{Wilson:isometries-for-NLG} any Riemannian manifold
$X$ admitting a transitive nilpotent Lie group of isometries is essentially
a nilpotent Lie group $N$ with an invariant metric. We shall determine
which six-dimensional (1-connected) nilpotent Lie groups $N$ generate
conformally parallel structures on manifolds of a special kind, described
hereby.

\section{Solvable extensions of nilpotent Lie algebras}
\label{sec:extension}

\noindent
Let $(N,h)$ denote a six-dimensional connected and simply-connected
nilpotent Lie group with a left-invariant Riemannian metric, and $\n$ its
Lie algebra. The orthonormal basis $\{ e^1,\dots , e^6 \}$ of the cotangent
bundle $T^*N$ is intended to be nilpotent, i.e.~such that
$de^i\in\Lambda^2V_{i-1}$, where the spaces $V_j =\text{span}_\mathbb
R\{e^1,\dots , e^{j-1}\}$ filtrate the dual Lie algebra: $0\subset
V_1\subset \ldots \subset V_5\subset V_6=\n^*$. The step-length of $\n$ is
defined as the number $p$ of non-zero subspaces appearing in the lower
central series\vspace{3pt}

\centerline{$\n\supseteq [\n,\n] \supseteq \bigl[ [\n,\n],\n\bigr]
\supseteq \ldots \supseteq \{0\}. $}\vspace{3pt}

\noindent
Given this, the terms Abelian and $1$-step are synonymous. We shall need
later the fact \cite{Nomizu:cohomology} that a nilmanifold
$\Gamma\backslash N$ and the Lie algebra of its universal cover have
isomorphic cohomology theories, $H^*(\n)\iso H^*_{\textrm{\tiny
    dR}}(\Gamma\backslash N)$.\\

\noindent Fix now a unit element $H\notin \n$ and suppose there exists a
non-singular self-adjoint derivation $D$ of $\n$ endowing
\begin{equation}
     \label{eq:extension}
     \s=\n\op \R H
\end{equation}
with the structure of a solvable Lie algebra. In other words think of
$\mathfrak s$ as an extension of the
following kind
\begin{definition}\label{Iwasawa extension}
     A metric solvable Lie algebra $\bigl( {\mathfrak s}, \lan \, , \, \ran
     \bigr)$ is said \emph{of Iwasawa type} if
     \begin{enumerate}
     \item $\mathfrak s= {\mathfrak a} \oplus {\mathfrak n}$, with
       $\n=[\mathfrak s,\mathfrak s]$ and
       $\mathfrak a=\n^\perp$ Abelian;
     \item $\ad_H$ is self-adjoint with respect to the scalar product $\lan\,
       , \,\ran$ and non-zero, for all $H\in\mathfrak a, H\neq 0$;
     \item for some (canonical) element $\tilde{H} \in\mathfrak a$,
       the restriction of $\ad_{\tilde H}$ to $\mathfrak n$ is
       positive-definite.\label{>0}
     \end{enumerate}
\end{definition}
\noindent
The terminology is clearly reminiscent of the Iwasawa decomposition of a
semisimple Lie group. This is indeed no coincidence, for any irreducible
symmetric space of non-compact type $Y=G/K$ can be isometrically identified
with the solvmanifold $S=AN$ relative to the decomposition $G=KAN$ of the
connected component of the isometry group of $Y$. Iwasawa-type extensions
are instances of standard solvmanifolds in the sense of Heber, and in a way
represent the basic model of standard Einstein
manifolds~\cite{Heber:noncompact-Einstein}. Now the nilpotent Lie groups of
concern (actually all, up to dimension six) always admit Einstein solvable
extensions~\cite{Lauret:standard-Einstein}, yet we wish to stress that all
known examples of non-compact homogeneous spaces with Einstein metrics are
of this kind, modulo isometries. What is more, they are completely
solvable, i.e.~the eigenvalues of any inner derivation are real. The
curvature of these spaces must be non-positive, because Ricci-flat homogeneous
manifolds are flat~\cite{Alekseevsky-K:Ricciflat-homogeneous}, and
Alekseevski\u{\i} has conjectured that a non-compact homogeneous Einstein
manifold has a transitive solvable isometry group. The latter cannot be
unimodular, as the space is assumed to be
non-flat~\cite{Dotti:unimodular-solvable}. This is in contrast to the
nilpotent picture, where a cocompact discrete subgroup always exists
\cite{Malcev:rational}, under the hypothesis of rationality of the
structure constants.

The whole point of reducing to rank one is that in the Einstein case, this is
no big specialisation, for~\cite{Heber:noncompact-Einstein} classifying
standard Einstein solvmanifolds is essentially the same as determining
those with $\codim [\mathfrak s,\mathfrak s]=1$.

Since $N$ has an invariant $\SU(3)$ structure one can suppose there exists
a diagonalisable operator $D\in \Lie{Der}(\n)$ with respect to a Hermitian
basis, that determines the rank-one extension as in~\eqref{eq:extension}.
That entails that there is indeed a unitary basis
consisting of eigenvectors --- let us still call it $\{e_i\}, i=1\ldots 6$
--- for which the matrix associated to $D=\ad_{\tilde H}$ is diagonal. Hence,
\begin{equation}
    \label{eq:D}
    \ad_{\tilde H} (e_i) = c_i e_i
\end{equation}
for some real constants $c_i$, which must be positive in order to satisfy
Definition~\ref{Iwasawa extension}. The derivation $D$ is chosen to
be precisely $\ad_{e_7}$, and
since the Cartan subalgebra $\mathfrak a$ is now
one-dimensional the only inner automorphism acting on $\n$ is the bracket
with the vector $\tilde{H}=e_7$, which is self-adjoint for the inner
product, and non-degenerate because $c_j\neq 0$, for all $j$'s.  Therefore,
the Maurer-Cartan equations of the rank-one solvable extension $\mathfrak
s=\n\oplus\R e_7$ assume the form
\begin{equation}
     \label{(*)}
     \left\{
       \begin{array}{l}
         de^j = \hd e^j + c_je^{j7}, \qq 1\leq j\leq 6\\
         d\7  = 0,
       \end{array}
     \right.
\end{equation}
where the `hat' indicates derivatives relative to the six-dimensional
world, i.e.~$\hd =d|_{\Lambda^*\R^6}$, and $\{e^j\}$ is the basis of $\s$
dual to $\{e_i\}$.  Results of Heber and Will
\cite{Heber:noncompact-Einstein,Will:Einstein-7-solv} guarantee that $Y=N
\times \R$ admits Einstein metrics, in fact there exists a unique choice of
the vector $(c_1,\ldots,c_6)$ such that the inner product $\lan \, , \,
\ran$ is Einstein.

In general, the Lie structure of $\n$ is defined by
\begin{equation}
     \label{eq:nilstructure}
     \left\{
       \begin{matrix}
         \hd\1  & = & a_1e^{12}    & + &  \ldots & + & a_{15}e^{56}\\
         \hd\2  & = & a_{16}e^{12} & + &  \ldots & + & a_{30}e^{56}\\
         \hdotsfor[2]{7}\\[2pt]
         \hd\6  & = & a_{76}e^{12} & + &  \ldots & + & a_{90}e^{56}
       \end{matrix}
     \right.
\end{equation}
\noindent
where all coefficients $a_{k}$ are real numbers.
\begin{rmk}\label{DJDJ=JDJD}
    This is a good point to see that the existence of a unitary basis
    diagonalising $\ad_{e_7}=D$ is tantamount to requiring that $D$ and $JDJ$
    commute, or $(DJ)^2=(JD)^2$.  This follows directly from
    \eqref{eq:standardstructure}, \eqref{eq:D}, and the computation for $e_1$
    is heuristic
    \begin{displaymath}
    \begin{array}{c}
     DJe_1=c_4e_4=\tfrac{c_4}{c_1}JDe_1, \q \textrm{hence} \\[3pt]
     (DJ)^2e_1=\tfrac{c_4}{c_1}DJ(JDe_1)=-\tfrac{c_4}{c_1}D^2e_1=
     -c_4c_1e_1= \hspace{5cm}\\
       \hspace{5cm} -\tfrac{c_1}{c_4}J(c_4^2e_4)=-\tfrac{c_1}{c_4}JD^2e_4=
       \tfrac{c_1}{c_4}JDDJe_1=(JD)^2e_1.
    \end{array}
    \end{displaymath}
    It has to be noticed though that $J$ need not necessarily be an almost
    complex structure for the argument. In fact, any endomorphism
    $\mathcal I$ of the tangent bundle of $N$ such that $\mathcal I\lan
    e_1\ran=\lan e_4\ran,\ \mathcal I\lan e_4\ran=\lan e_1\ran$ {\sl et
      cetera} does the job, since then $\mathcal I^{-1}D\mathcal I$ and $D$
    are simultaneously diagonalisable, hence commute.
   \end{rmk}

\section{The classification}\label{sec:classification}

\noindent
Before we start investigating equations \eqref{eq:m} in relation to the
induced geometry on $N^6$, let us discuss the delicate point of the choice
of the $\SU(3)$ reduction.   Define $\psi^{\pm}$ as in
\eqref{eq:standardstructure} with $\{ e_i \}$  a unitary basis that 
diagonalises $D$. A  reduction to $\SU(3)$ is
determined by the choice of an element
$$
  {\tilde\psi}^+  = \psp\cos\theta + \psm\sin\theta
$$
(for some angle $\theta$)
in the circle generated by $\psi^+$ and $\psi^ -$ in $\W_1$.  In 
general, it is impossible to express $\tilde \psi^+$  in
terms of a basis that diagonalises $D$ as simply as in 
$\eqref{eq:standardstructure}$.  The proof of next theorem shows
  that one can   in fact  assume that $\theta = 0$.

Moreover,  one can say is that there is a unique --- up to sign --- closed
$3$-form in the circle  (Proposition \ref{Jintegrable}).

Let us write
\begin{displaymath}
     \f=\om\7+\psp, \qquad \ff=\psm\7+\tfrac 12\oo,
\end{displaymath}
whence one immediately finds that
\begin{equation}
     \label{(**)}
     d\om\7+d\psp=-3m\psp\7,\qquad d\psm\7+\omega d\omega =2m\oo\7.
\end{equation}
Reflecting the splitting of the fibres of the cotangent bundle
$T_y^*Y=\R^6\oplus \R\7$, the relations \eqref{(*)} give
\begin{displaymath}
     d\om=\hd\om-\bigl((c_1+c_4)e^{14}-(c_3+c_2)e^{23}+(c_5+c_6)e^{56}\bigr)\7,
\end{displaymath}
so
\begin{displaymath}
    \begin{split}
     d\oo & = \hd\oo + 2\bigl((c_1+c_4+c_3+c_2)e^{1423} +
     (c_3+c_2+c_5+c_6)e^{2356} \\
            & \phantom{MMMMMM} - (c_1+c_4+c_5+c_6)e^{1456}\bigr)\7.
    \end{split}
\end{displaymath}
Similarly one computes the exterior derivatives of the real $3$-forms:
\begin{displaymath}
     \begin{split}
       d\psp = \hd\psp & +(c_1+c_2+c_5)e^{1257}+(c_1+c_3+c_6)e^{1367}\\
                       & -(c_3+c_4+c_5)e^{3457}+(c_2+c_4+c_6)e^{2467},\\
       d\psm = \hd\psm & +(c_1+c_2+c_6)e^{1267}-(c_1+c_3+c_5)e^{1357}\\
                       & -(c_2+c_4+c_5)e^{2457}-(c_3+c_4+c_6)e^{3467}.
     \end{split}
\end{displaymath}

\bigbreak\noindent
When, in general, $G_2$-manifolds $Y$ are constructed starting from six
dimensions, many of their features are determined by the underlying
$\SU(3)$ structure, and the following definition becomes natural
\begin{definition}
    {\rm\cite{Chiossi-S:SUG}} An almost Hermitian manifold is
     \emph{half-flat}, or \emph{half-integrable}, if the reduction
     is such that both $\psp$ and $\oo$ are closed (with respect to $\hd$).
\end{definition}
\noindent
This is the same as asking that the intrinsic torsion components
$W_1^+,W_2^+,W_4$ and $W_5$ vanish simultaneously. This sort of structure
appears, in various disguises, on any hypersurface in $\R^7$ (or Joyce
manifold, for that matter), and its possible r\^ole in $\mathcal M$-theory
has been recently examined~\cite{Gurrieri-LMW:mirrorsymmetry,
    Becker-DKT:duality-cycle}.

\smallbreak
Plugging the previous equations into system
\eqref{eq:m} allows us to discover a geometrical constraint, for
\begin{lemma}
    When $(Y, \varphi)$ is conformal to a $G_2$-holonomy manifold, $N$ has a
    half-flat $\SU(3)$ structure.
\end{lemma}
\begin{proof}
This is clear if one considers the terms in \eqref{(**)}
that belong to $(\7)^\perp$.
\end{proof}

\smallbreak \noindent
On the other hand, the components of \eqref{eq:m} in the
direction of $\7$ read
\begin{equation}
     \label{eq:m-equations}
     \left\{
       \begin{array}{lcl}
        \hd\om  &\!\! =\!\! & -(c_1+c_2+c_5+3m)e^{125}-
                              (c_1+c_3+c_6+3m)e^{136}- \\[3pt]
                &           & \phantom{-}(c_3+c_4+c_5+3m)e^{345}-
                              (c_2+c_4+c_6+3m)e^{246} \\[3pt]
        \hd\psm &\!\! =\!\! & (2m-c_1-c_2-c_3-c_4)e^{1423} +
                              (2m-c_3-c_2-c_5-c_6)e^{2356}+\\[3pt]
                &           & (2m+c_1+c_4+c_5+c_6)e^{1456}.
       \end{array}
     \right.
\end{equation}

\smallbreak\noindent
We will show that the derivation $D = ad_{e_7}$ has an eigenvector (for
instance $e_1$) belonging to $[{\mathfrak n}, {\mathfrak n}]^{\perp}$, so
the structure of $\s$ is determined by equations \eqref{(*)}, with $\hd e^1
=0$ and $\hd e^j$ given by \eqref{eq:nilstructure} for $j=2,\ldots,6$. The
point is to find all possible coefficients $a_{k}, k=16,\ldots, 90$ and
$c_j, 1\leq j\leq 6$ such that $d^2(e^j)=0$ and \eqref{(**)} are satisfied,
for some non-vanishing $m$. In this way we obtain the following classifying
result
\begin{theorem}
    \label{theorem}
    Let $N$ be a nilpotent Lie group of dimension $6$ endowed with an
    invariant $\SU (3)$ structure $(\omega, {\tilde \psi}^+)$. Suppose
    there is a non-singular and self-adjoint derivation $D$ of the Lie
    algebra $\n$ such that $(DJ)^2 = (JD)^ 2$. Then on the solvable
    extension $\s=\n\oplus \R e_7$ with $\ad_{e_7}=D$, the $G_2$ structure
    $\varphi=\omega\wedge\7+ {\tilde\psi}^+$ is conformally parallel if and
    only if $\n$ is isomorphic to one of the following:
    \begin{center}
      \begin{tabular}{lr}
       $(0,0,0,e^{12},e^{13},e^{23})$, & $(0,0,0,0,e^{12},e^{13})$,\\[4pt]
       $(0,0,0,0,e^{12},e^{14}+e^{23})$, & $(0,0,0,0,0,e^{12}+e^{34})$,\\[4pt]
       $(0,0,0,0,e^{13}+e^{42},e^{12}+e^{34})$, & $(0,0,0,0,0,e^{12})$,\\[4pt]
       $(0,0,0,0,0,0)$.
      \end{tabular}
    \end{center}
  \end{theorem}
\noindent
Though the list does not appear that meaningful at first sight, it becomes
more significant once considered in relation to the descriptions given in
\cite{Magnin:nil-7, Goze-K:nilpotent}.
\begin{proof}
  Since $D$ is a derivation, it must preserve the orthogonal splitting
  $[\n,\n]\op [\n,\n]^\perp$.  Indicating the derived algebra $[\n,\n]$ by
  $\n^1$, one infers that $\ad_{e_7}(\n^1) \set \n^1$, hence
  $\ad_{e_7}(\n^1)^{\perp} \set (\n^1)^{\perp}$. Then one can suppose that
  there exists a unitary basis $\{e^i\}$ which diagonalises $D$ with $\1$
  closed in $\n^*$.  The structure equations of $\s$ are given by
  \eqref{(*)} and \eqref{eq:nilstructure}, with $a_j=0$ for all
  $j=1,\ldots,15$ and the $3$-form $\tilde \psi^+$ can be then expressed,
  in terms of the previous basis $\{e^i\}$, as
$$
\psi^+ \cos \theta + \psi^ - \sin \theta,
$$
for some angle $\theta$, where $\psi^{\pm}$ are given by
$\eqref{eq:standardstructure}$.  It is necessary to impose the quadratic
relations $d^2e^i=0$ together with the linear equations \eqref{eq:m} , for
a total of $35\cdot 5+56$ constraints. The complete system has $75+6+1+ 1$
unknown variables $a_k, c_j, m, \theta$.  Given the number of parameters
and equations, the results were also checked with the {\sc Maple} package.
Inserting the coordinates, \eqref{eq:m} yields a bulk of $56$ linear
constraints on $83$ coefficients. We have to distinguish two cases: $\theta
= 0$ and $\theta \neq 0$.  If $\theta$ is not zero, the imposition of
$d^2e^i=0$ and $c_l \neq 0$ gives no solution.

If $\theta = 0$ the equations  \eqref{eq:m} reduce to
   \eqref{eq:m-equations}. Replacing the respective
     expressions in~\eqref{eq:nilstructure}, one gets
     {\allowdisplaybreaks\begin{align*}
       de^1  &= c_{1}e^{17} ,\\
       de^2  &= c_{2} e^{27} +a_{16} e^{12}+ (a_{65} -a_{22}
       -a_{89}) e^{13} +
         a_{55} e^{14}+(a_{77}+a_{56})e^{15} + a_{20}e^{16}+\\
         & b_1 e^{23} + a_{22}e^{24}+ a_{23}e^{25} + a_{24}e^{26}
         +(a_{88}-a_{32}-a_{48}-a_{64})e^{34} + a_{26}e^{35}+\\
         & a_{27}e^{36} +
         (2a_{62}+2c_1-2c_5+a_{53}+a_{20}+a_{67}-a_{76})e^{45}
         + a_{29}e^{46}+ a_{30}e^{56},\\
         de^3  &= c_{3}e^{37}+(-a_{22}+a_{89}-a_{65})e^{12}+
         a_{32}e^{13}+ a_{33}e^{14}+\\
         & (4c_1-2c_2-2c_5+2a_{53}+3a_{62}+3a_{20}-a_{76})e^{15} +
         a_{35}e^{16} + a_{36}e^{23}+\\
         & (-a_{88}+a_{16}+a_{48}+a_{64})e^{24}+ (-a_{26}+2a_{24}+
         a_{66})e^{25}+ a_{39}e^{26} -a_{22}e^{34}+\\
         & b_2e^{35}+ b_3e^{36}+(a_{29}-a_{77}-a_{56}-a_{35})e^{45}+\\
         & (5a_{20}+6c_1-4c_2-4c_5+4a_{53}+6a_{62}-a_{76}+2c_3+
         a_{85}-a_{57})e^{46}+ a_{45}e^{56},\\
         de^4  &= c_4 e^{47} +a_{55}e^{12}+ a_{33}e^{13}+ a_{48}e^{14}+
         a_{49}e^{15} + a_{50}e^{16}+ \\
         & (-a_{16}-a_{32})e^{23} -a_{33}e^{24}
         + a_{53}e^{25}+ a_{56}e^{26} + a_{55}e^{34}+\\
         & a_{56}e^{35} + a_{57}e^{36} -a_{50}e^{45} + a_{49}e^{46}
         + (-a_{88}+2a_{64}-a_{74})e^{56},\\
         de^5  &= c_{5}e^{57} +(a_{35}+a_{56})e^{12}+ a_{62}e^{13}+
         a_{49}e^{14}
         +a_{64}e^{15} + a_{65}e^{16} + a_{66}e^{23}  + a_{67}e^{24}+\\
         & b_4 e^{25}  + (-a_{83}+2a_{71}-a_{30})e^{26}
         + a_{29}e^{34}   + a_{71}e^{35}+ a_{72}e^{36} - a_{89}e^{45}+\\
         & a_{74}e^{46}+ b_5 e^{56},\\
         de^6  &= c_6 e^{67}+a_{76}e^{12} + a_{77}e^{13}+ a_{50}e^{14}+
         (2a_{89}-a_{65})e^{15} + b_6 e^{16}+\\
         &    (2a_{23}+a_{27}+a_{39})e^{23} + (a_{29}-a_{77}-a_{56}-
         a_{35})e^{24} + a_{83}e^{25} + b_7 e^{26}+\\
         &    a_{85}e^{34}+ (2a_{33}-2a_{36}+a_{72}+a_{45})e^{35} +
         b_8  e^{36} + a_{88}e^{45} + a_{89}e^{46}+b_9 e^{56},\\
         de^7   &=  0,
       \end{align*}}
       in terms of a certain number of parameters. We have put
         \begin{gather*}
           b_1= -a_{83}+a_{55}+a_{71}, \qquad b_2=a_{27}+a_{39}+a_{23}, \qquad
           b_3=-a_{24}-a_{66}, \\
           b_4=-a_{33}-a_{72}+a_{36}-a_{45}, \quad 
           b_5=-a_{49}+a_{24}-a_{26}+a_{66},
           \quad b_6=-a_{88}+a_{64}-a_{74}, \\
           b_7=a_{33}+a_{72}-a_{36}, \qquad b_8=-a_{71}+a_{30}, \qquad
           b_9 = a_{39}+a_{23}-a_{50}
         \end{gather*}
         for convenience. Besides, the following relations must hold:
         \begin{gather*}
           c_4= -5c_1+3c_2+4c_5-4a_{53}-5a_{62}-4a_{20}-c_3-a_{67}+a_{76}-a_{85},\\
           c_6 = 3c_{2}-c_3+3c_5+a_{57}-3a_{53}-4c_1-4a_{62}-4a_{20},\\
           m = c_1 - c_2 - c_5 + a_{53} + a_{62} + a_{20}.
         \end{gather*}
         Only at this point it seems realistic to annihilate the quadratic
         relations coming from the Jacobi identity, hence set to zero the
         coefficients of the terms $e^{ij7}$ appearing in the various $d^2e^i=0$.
         Since $c_j \neq 0$ for all $j=1,\ldots,6$, the closure of
         $de^i$ kills all $b_l$'s above, and furthermore\\
         \centerline{$ a_{16} = a_{22} = a_{24} = a_{23} = a_{32} =
           a_{33} = a_{36} = a_{48} = a_{49} = a_{50} = a_{55} = a_{64} =
           a_{71} = a_{89} = 0.  $}
         Thus, the structure equations eventually
         reduce to a simpler form
         \begin{displaymath}
           \begin{array} {lcl}
             de^1 & = & c_{1}e^{17} \\[3pt]
             de^2 & = & c_{2}e^{27}
             +a_{65}e^{13}+(a_{77}+a_{56})e^{15}+a_{20}e^{16}-a_{74}e^{34}+\\
             & & (2a_{62}+2c_1-2c_5+a_{53}+a_{20}+a_{67}-a_{76})e^{45}+
             a_{29}e^{46}\\[3pt]
             de^3 & = & c_{3}e^{37}
             -a_{65}e^{12} +(4c_1-2c_2-2c_5+2a_{53}+
             3a_{62}+3a_{20}-a_{76})e^{15}+\\[3pt]
             & &  a_{35}e^{16} + a_{74}e^{24}  +
             (a_{29}-a_{77}-a_{56}-a_{35})e^{45} +\\[3pt]
             & &  (5a_{20}+6c_1-4c_2-4c_5+4a_{53}+6a_{62}-a_{76}+
             2c_3+a_{85}-a_{57})e^{46}\\ [3pt]
             de^4 & = & c_4 e^{47}
             +a_{53}e^{25}+ a_{56}e^{26} + a_{56}e^{35} +a_{57}e^{36}\\[3pt]
             de^5 & = & c_{5}e^{57}
             +(a_{35}+a_{56})e^{12}+ a_{62}e^{13}+
             a_{65}e^{16} + a_{67}e^{24}+ a_{29}e^{34} + a_{74}e^{46}\\[3pt]
             de^6 & = & c_6 e^{67}
             +a_{76}e^{12} + a_{77}e^{13} -a_{65} e^{15}+
             (a_{29}-a_{77}-a_{56}-a_{35})e^{24}+ \\[3pt]
             & &  a_{85}e^{34} - a_{74}e^{45} \\ [3pt]
             de^7 & = & 0.
           \end{array}
         \end{displaymath}
         In particular, the vanishing of the
         coefficients of $e^{137}, e^{347}$ in $d^2 (e^2)$ and of $e^{127}, e^{247}$
         in $d^2 (e^3)$ yields
         \begin{displaymath}
           \begin{array} {l}
             a_{65} (c_1 - c_2 + c_3) = 0, \quad a_{74} (- c_2 - c_3 + c_4) = 0,\\
             a_{65} (- c_1 - c_2 + c_3) = 0, \quad a_{74} ( c_2 - c_3 + c_4) = 0,\\
           \end{array}
         \end{displaymath}
         thus $a_{65} = a_{74} = 0$. The terms $e^{357}, e^{267}$ in $d^2 (e^4)$
         similarly give
         \begin{displaymath}
           \begin{array} {c}
             a_{56} (- c_3 + c_4 - c_5) = 0, \quad a_{56} (- c_2 + c_4 - c_6)=0,
             \quad a_{53}  ( - c_2 + c_4 - c_5) = 0,\\
             a_{57} (- c_2 + a_{57} + a_{53} + c_1 + a_{62} + c_3 + a_{67} -
             a_{76} + a_{85}) = 0.
           \end{array}
         \end{displaymath}
         Altogether, the following cases crop up. We shall examine them
         one by one trying to make further coefficients disappear.
         
         \bigbreak
         \noindent
         {\bf Case a)} $a_{53} = a_{56} = a_{57} = 0$ (corresponding to
         $\hd e^4 = 0$).
         The relation $d^2=0$ yields $a_{29} = a_{35} = a_{77} =0$, and six
         non-isomorphic algebra types come out:
         \begin{equation} \label{abelian}
           (-me^{17},-me^{27},-me^{37},-me^{47},-me^{57},-me^{67},0)
         \end{equation}
     has an underlying Abelian Lie algebra, if one disregards the
     $D$-action. Next,
     \begin{equation}
       \label{h3}
       (-\tfrac 23 me^{17},-me^{27},-\tfrac 43 me^{37}+\tfrac 23
       me^{15}, -me^{47}, -\tfrac 23 me^{57},-me^{67},0)
     \end{equation}
     extends $\n \iso (0,0,0,0,0, e^{12})$;
     \begin{equation}
       \label{12+34}
       \bigl(-\tfrac 34 me^{17},-me^{27},
       -\tfrac 32 me^{37}+\tfrac 12 m(e^{15}-e^{46}),-\tfrac 34 me^{47},
       -\tfrac 34 me^{57},-\tfrac 34 me^{67},0\bigr)
     \end{equation}
     is clearly given by $\n \iso (0,0,0,0,0,e^{12}+e^{34})$;
     \begin{equation}
       \label{12,13+24}
       \bigl(-\tfrac 45 me^{17},-\tfrac 65 me^{27}-\tfrac 25 me^{45},
       -\tfrac 75 me^{37}+\tfrac 25 m(e^{15}-e^{46}),-\tfrac 35 me^{47},
       -\tfrac 35 me^{57},-\tfrac 45 me^{67},0\bigr)
     \end{equation}
     attached to $\n \iso (0,0,0,0,e^{12},e^{14}+e^{23})$;
     \begin{equation}
       \label{12,13}
       (-m e^{17},-\tfrac 54 me^{27}-\tfrac 12 me^{45},
       -\tfrac 54 me^{37}-\tfrac 12 me^{46}, -\tfrac 12 me^{47},-\tfrac
       34 me^{57}, -\tfrac 34 me^{67},0),
     \end{equation}
     whose $\n$ is essentially $(0,0,0,0,e^{12},e^{13})$;
     \begin{equation}
       \label{13+42,12+34}
       \bigl(-\tfrac 23 me^{17},-\tfrac 43 me^{27}-\tfrac 13 m(e^{16}+e^{45}),
       -\tfrac 43 me^{37}+\tfrac 13 m(e^{15}-e^{46}),-\tfrac 23 me^{47},
       -\tfrac 23 me^{57},-\tfrac 23 me^{67},0\bigr)
     \end{equation}
     is an extension of the Iwasawa Lie algebra, isomorphic to
     $(0,0,0,0,e^{13}+e^{42},e^{12}+e^{34})$.

     \bigbreak
     \noindent
     {\bf b)} $a_{53} = a_{56} = 0, a_{57} = c_{5} - c_1 - a_{62} - c_3 -
     a_{67}+ a_{76} - a_{85}$. Up to isomorphism, one gets the two Lie algebras
     with structure~\eqref{h3} and~\eqref{12,13}.

     \bigbreak
     \noindent
     {\bf c)} $a_{56} =  a_{57} = 0, c_2 = \tfrac 52
     c_1- \tfrac 32 c_5 +\tfrac 52 a_{62} + 2 a_{20} +\tfrac 12
     c_3 + \tfrac 12 a_{67}+2 a_{53} - \tfrac 12 a_{76} +
     \tfrac{1} {2} a_{85}$. This time around one finds the algebras of b) plus
     \begin{equation}
       \label{12,13,23}
       (-\tfrac 35 me^{17},-\tfrac 35 me^{27},-\tfrac 65 me^{37}+\tfrac
       25 me^{15},
       -\tfrac 65 me^{47}+\tfrac 25 me^{25},-\tfrac 35 me^{57},
       -\tfrac 65 me^{67}+\tfrac 25 me^{12},0),
     \end{equation}
     which arises from $\n \iso (0,0,0,e^{12},e^{13},e^{23})$.
     \bigbreak
     \noindent
     {\bf d)} $a_{56} = 0$, $c_2 = \tfrac 52 c_1 -\tfrac 32 c_5 + \tfrac 52
     a_{62} +2 a_{20} + \tfrac 12 c_3 + \tfrac 12 a_{67} + 2 a_{53} - \tfrac
     12 a_{76} + \tfrac 12 a_{85}$, $a_{57} = c_5 -a_{53} -c_1-a_{62} - c_3
     -a_{67} +a_{76} -a_{85}$, by which one regains~\eqref{12+34}.

     \bigbreak
     \noindent
     {\bf e)} $a_{53} = a_{57} = 0$,
     $a_{20} =-c_1+c_2+\tfrac 12 c_5 -a_{62}-\tfrac 12 c_3$,
     $a_{76} = c_1+c_2-c_5+a_{62}+a_{67}+a_{85}$
     yield no solutions.

     \bigbreak
     \noindent
     {\bf f)} $a_{76} = c_1+c_2-c_5+a_{53}+a_{57}+a_{62}+a_{67}+a_{85}$, $c_3 =
     c_2$, $a_{20} = -c_1+ \tfrac 12 c_2 + \tfrac 12 c_5 - \tfrac 34 a_{53}
     -a_{62} + \tfrac 14 a_{57}$. The last case produces~\eqref{12+34} one more
     time, and basically concludes the proof of the Theorem.
\end{proof}

\medbreak\noindent
The fact that so few Lie algebras are gotten may depend on the requirements
made both on the $G_2$ structure and on the seven-dimensional construction.
One easily recognizes that the groups associated to \eqref{abelian} and
\eqref{h3} are the torus $T^6$ and the product $T^3\times
\mathcal H^3$ of a torus with the real $3$-dimensional Heisenberg group
respectively, while \eqref{13+42,12+34} is attached to the complexified
$3$-dimensional Heisenberg group ${\mathcal H}_3^\C$.

In the case of an Einstein solvmanifold, the eigenvalues of the derivation
$D$ are positive integers $k_1<\ldots <k_r$ without common
divisors~\cite{Heber:noncompact-Einstein}. Indicating the respective
multiplicities by $d_1, \dots, d_r$, Heber defines the string $(k_1,
\ldots, k_r; d_1, \ldots, d_r)$ the eigenvalue type of the solvmanifold.
In our situation we have the following

\begin{corollary} 
With the above hypotheses, the eigenvalues  of the 
  derivation $D$ have all the same sign. More precisely, the eigenvalue
  data of the algebras of Theorem \ref{theorem} are given by the following
  scheme

\bigbreak
\begin{center}
    \begin{tabular}{|c|c|c|}
      \hline
      $\vs \hbox{Nilpotent Lie algebra}$ & \hbox{Eigenvalues} &
      \hbox{Multiplicities}\\
      \hline\hline
      $\vs (0,0,e^{15},0,0,0)$ &$ -\tfrac 23 m, -m, -\tfrac 43 m$ & $2,3,1$\\
      \hline
      $\vs (0,0,e^{15} + e^{64},0,0,0)$ &
      $-\tfrac 34 m, -m, -\tfrac 32 m$ & $4,1,1$\\
      \hline
      $\vs (0,e^{45},e^{64}+e^{51},0,0,0)$ &
      $-\tfrac 35 m,-\tfrac 45 m,-\tfrac 65 m,-\tfrac 75 m$ & $2,2,1,1$\\
      \hline
      $\vs (0,e^{45},e^{46},0,0,0)$ &
      $-\tfrac 12 m, -\tfrac 34 m, -m, -\tfrac 54 m$ & $1,2,1,2$\\
      \hline
      $\vs (0,e^{16}+e^{45},e^{15}+e^{64},0,0,0)$ &
      $-\tfrac 23 m, -\tfrac 43 m$ & $4,2$\\
      \hline
      $\vs (0,0,e^{15},e^{25},0,e^{12})$ & $-\tfrac 35 m,-\tfrac 65 m$ 
& $3,3$ \\
      \hline
    \end{tabular}
\end{center}\vspace{0.4cm}
\end{corollary}

The real number  $m$ has to be
a negative, in order for $\s$ to be of Iwasawa type. Notice
that the result holds just assuming non-degeneracy (thus dropping
(\ref{>0}) in Definition \ref{Iwasawa extension}).

In the present set-up, the Corollary matches to
the result of \cite{Will:Einstein-7-solv} for appropriate choices of $m$.
Moreover, the eigenvalue type is unique to each example, in contrast to the
Einstein case where the solvmanifolds associated $\h_3\oplus\h_3$ and 
$\h_3^\C$ have the same eigenvalue
type.

\medbreak
As a by-product of the classification, the following necessary 
condition crops up:
\begin{corollary}
With the above hypotheses, if $Y$ has a $G_2$ structure of type $\X_4$,
then $N$ is either a 2-step nilpotent Lie group or a torus.\qed
\end{corollary}
\noindent
This is reflected in the fact that the metrics supported by these
solvmanifolds arise on torus bundles over tori of various dimensions and
rank, cf.~\S \ref{sec:metrics}. Notice that the only 2-step nilmanifold
missing, so to speak, is that corresponding to
$(0,0,0,0,e^{12},e^{34})\iso\h_3 \oplus \h_3$.  It is known to the authors
that this Lie algebra admits a large family of half-flat
$\SU(3)$-structures. The product of the corresponding nilmanifold with some
real interval can be endowed with a metric with holonomy contained in $G_2$
\cite{Hitchin:forms}, but by the Theorem such metric will not be
conformally equivalent to a homogeneous one on a solvable extension of
${\mathcal H}^3 \times {\mathcal H}^3$.

\subsection{Some consequences}

\noindent
A slight change of approach allows to detect properties in a simpler way.
Let
\begin{displaymath}
     \alpha^1=\1+i\4,\quad \alpha^2=\2-i\3,\quad \alpha^3=\5+i\6
\end{displaymath}
be the basis of complex $(1,0)$-forms determined
by~\eqref{eq:standardstructure}, so one may write
\begin{displaymath}
     \Psi=\alpha^1\wedge\alpha^2\wedge\alpha^3=\alpha^{123}.
\end{displaymath}
Translate all $\SU(3)$ forms into this language
\begin{displaymath}
    \om=-\tfrac 1{2i}(\alpha^{1\bar 1}+\alpha^{2\bar 2}+\alpha^{3\bar 3})
    \q \text{and} \q
    \psp=\tfrac 12 (\alpha^{123}+\alpha^{\bar 1\bar 2\bar3}) \q \textsl{etc.},
\end{displaymath}
with $\alpha^{\bar{\imath}}$ indicating the conjugate
$\overline{\alpha^\imath}$ and $\alpha^{ij}$ standing for
$\alpha^i\alpha^j,\ i,j=1,2,3$. Equations~\eqref{eq:m-equations} become
\begin{gather*}
    \begin{array}{rcl}
      2i\hd(\alpha^{1\bar 1}+\alpha^{2\bar 2}+\alpha^{3\bar 3})
      & = & (c_4-c_1)(\alpha^{123}+\alpha^{\bar 1\bar 2\bar 3})+
       (c_3-c_2)(\alpha^{12\bar 3}+\alpha^{\bar 1\bar 2 3})+ \\[3pt]
      &   & (c_6-c_5)(\alpha^{1\bar 2 3}+\alpha^{\bar 1 2\bar 3})\ -
            (6m+\textstyle\sum_i c_i)(\alpha^{\bar 1 23}+\alpha^{1\bar 
2\bar 3}),
    \end{array}\\[5pt]
    \begin{array}{rcl}
       -4i\hd\alpha^{123}
       & = & (2m-c_1-c_2-c_3-c_4)\alpha^{1\bar 12\bar 2}\ +
             (2m-c_3-c_2-c_5-c_6)\alpha^{2\bar 23\bar 3}\ - \\[3pt]
       &   & (2m+c_1+c_4+c_5+c_6)\alpha^{1\bar 13\bar 3}.
    \end{array}
\end{gather*}
Because of type, the terms in the latter can be treated separately
\begin{equation}\label{d-alphas}
     \left\{
       \begin{array}{rcl}
         4id\alpha^3 & = & (2m-c_1-c_2-c_3-c_4)\alpha^{\bar 1\bar 2}\\
         4id\alpha^2 & = & (2m+c_1+c_4+c_5+c_6)\alpha^{\bar 1\bar 3}\\
         4id\alpha^1 & = & (2m-c_2-c_3-c_5-c_6)\alpha^{\bar 2\bar 3}.
       \end{array}
     \right.
\end{equation}
The special case in which $J$ is actually a complex structure is
instructive. Since the six-dimensional manifold is half-flat, hence has
intrinsic torsion only in $\W_1^-, \W_2^-$ and $\W_3$, the further
requirement that $\hd\psm=0$ forces $N$ to become balanced.
By~\eqref{eq:m-equations} the extension's coefficients satisfy the relation

\smallbreak
\centerline{$ c_1+c_4 = c_5+c_6 = -\tfrac 13 (c_2+c_3).$}

\smallbreak
\noindent
In relation to the structures \eqref{abelian} -- \eqref{12,13,23},
equations~\eqref{eq:m} confirm that if $(N,J)$ is Hermitian then $Hol(g)$
is a subgroup of $G\sb 2$
\begin{proposition}
   \label{Jintegrable}
    On the solvable group corresponding to $\mathfrak s=\n\oplus \R e_7$ the
    $G\sb 2$ structure $\f=\omega\7+\psp$ cannot be conformally parallel if
    the almost complex structure $J$ on $\n$ is integrable.\qed
\end{proposition}
By this result, if one fixes the almost complex structure
$J$, there exists a unique choice of $\psi^+$ in the $\U(1)$-family
of stable real $3$-forms, since the closure of $\psi^-$
renders $J$ necessarily integrable.

\medbreak
\begin{rmk}
    A reasonable question is to ask whether this geometry
    has the potential to produce strong $G_2$-metrics~\cite{Cleyton-S:sG2}.
    As the torsion form is merely $\Phi=m\psm$, its closure entails that $N$ is
    again a complex manifold (hence balanced). Although this is enough to
    conclude that the holonomy of $Y$ reduces, things get even worse, for
    $d\psm$ has components in $\Lambda^3\n^* \wedge \7$ as well, forcing

     \smallbreak
     \centerline{$ c_1=c_4, c_5=c_6, c_2=c_3=-(c_1+c_5).$}

     \smallbreak
     \noindent
     This gives $\hd\omega=-3m\psp \in \W_1$, whence $N$ has to be K\"ahler,
     confirming that the only solutions come from taking $m=0$ in
     \eqref{eq:m}. We conclude that if $dT=0$, the Lie algebra structure
     \eqref{(*)} simplifies to $de^j=\hd e^j, d\7=0$, so we are merely
     looking at $Y$ as the Riemannian product of $N$ with $\R$, much of which
     is known~\cite{Chiossi-S:SUG}.
\end{rmk}

\bigbreak
A similar argument also restricts the range
of $m$ in the general set-up. We claim in fact that\vspace{4pt}

\centerline{
     $ 2m \in \{ c_1+c_2+c_3+c_4,\; -(c_1+c_4+c_5+c_6), \;
     c_2+c_3+c_5+c_6\}. $
}\vspace{4pt}

\noindent
If $\s$ does not satisfy the above relation, then all coefficients in
\eqref{d-alphas} are different from zero, affecting the topology of $N$.
Considering $\alpha^2$ for instance, one sees that $\2$ and $\3$ cannot be
simultaneously closed, so the manifold $N$ cannot admit more than three
independent closed 1-forms. But Theorem~\ref{theorem} tells that only Lie
algebras with first Betti number $b_1\geq 3$ crop up, and the unique
$2$-step algebra attaining the minimum is \eqref{12,13,23}, which fails to
satisfy the assumption.

\medbreak The conditions to have a compatible almost K\"ahler structure are
found in a similar fashion. It is non-obvious, and certainly unusual, that
the symplectic condition also annihilates the component of the intrinsic
torsion in $\W_2^-$, a module not directly depending upon $d\omega$:
\begin{proposition}
    Under the above assumptions, the nilpotent Lie group $(N,J,\omega)$ is
    symplectic only when it is a torus, in other words

\centerline{$d\omega=0 \iff \n\ {\rm is\ Abelian.}$}
\end{proposition}
\begin{proof}
   If $\om$ is closed, its expression in complex form easily gives $c_1=c_4,
   c_3=c_2, c_5=c_6$, which corresponds precisely to
   $J\ad_{e_7}=\ad_{e_7}J$; concerning Remark~\ref{DJDJ=JDJD}, it is
   definitely worth noticing that the almost complex structure and the
   derivation $D=\ad_{e_7}$ commute just for two nilpotent Lie algebras,
   that is the Abelian one and the Iwasawa Lie algebra. The latter though
   does not satisfy the requirement that $m = -\sum c_i = -\tr\,\ad_{e_7}$,
   whence only the torus $T^6$ has symplectic structures generating a
   $G_2$-manifold of type $\X_4$.
\end{proof}

\section{Description of the Ricci-flat metrics}
\label{sec:metrics}

\noindent
So $\mathfrak s = \mathfrak n \oplus \R e_7$ possesses a conformally parallel
$G_2$ structure determined by the Lie types~\eqref{abelian} --
\eqref{12,13,23} of $\n$. A transformation $g\tto e^{2f}g$ with conformal
factor $d f = -m e^7$ produces Ricci-flat metrics,
which we describe in detail. Since the corresponding simply-connected
solvable Lie group Lie group $S$ is diffeomorphic to $\R^7$, it is possible
to find global coordinates $(x_1, \ldots x_6, t)$ that describe the
left-invariant $1$-forms $\1,\ldots,\6$ and $\7=dt$, all of which depend
upon one real parameter $m\neq 0$. The general form for these metrics will
thus be $ g = e^{-2 m t} \sum_{i = 1}^7 (e^i)^2$, and the explicit
calculations will be relevant in the determination of holonomy groups.

In some cases, the solution can be related to the results of
\cite{Gibbons-LPS:domain-walls}, whose metrics depend upon a function
accounting for the scaling symmetry. We can thus prove that all our metrics
admit a homothetic Killing field, i.e.~a vector field $Z$ such that
${\mathcal L}_Z g= cg, c\in\R$.

\subsection{The Abelian case}

For the Lie algebra~\eqref{abelian}, the coordinates
\begin{displaymath}
     \left\{
       \begin{array} {l}
         e^i = e^{m t} dx_i,\quad i = 1, \ldots, 6,\\
         e^7 = dt
       \end{array}
     \right.
\end{displaymath}
just yield the flat metric $ g = \sum_{i = 1}^6 dx_i^2 + e^{-2 m t} dt^2 $
on $T^6\times\R$.

\medbreak
\subsection{The algebra $\boldsymbol{\R^3 \oplus \h_3}$}

Consider the solvable extension of the product~\eqref{h3} of a torus $T^3$
with ${\mathcal H}^3$, and let
\begin{equation}
     \label{lambda coefficients}
     \left\{
       \begin{array} {l}
         e^i = e^{\tfrac 23 m t} dx_i,\quad i = 1, 5,\\
         e^l = e^{m t} dx_l,\quad l = 2, 4, 6,\\
         e^3 = -\tfrac 23 m\, e^{\tfrac 43 m t} ( dx_3 + x_5 dx_1),\\
         e^7 = dt.
       \end{array}
     \right.
\end{equation}
The Riemannian structure
\begin{equation}
     \label{7metric}
       g = dx_2^2 + dx_4^2 + dx_6^2
       + e^{-\tfrac 23 m t} ( dx_1^2 + dx_5^2)
       + \tfrac 49 m^2 e^{\tfrac 23 m t} (dx_3 +  x_5 dx_1)^2 + e^{-2m t} dt^2
\end{equation}
restricts to a special holonomy metric $ds$ on
$\textsl{span}\{x_1,x_3,x_5,t\}$ viewed as $Q \times \R$, $Q$ being the
total space of a circle bundle over $T^2$. In fact the subgroup
of $G_2$ preserving $ds$ is orthogonal, hence $Hol(g)=G_2\cap
\SO(4)=\SU(2)$.

\medbreak
\subsection{The algebra $\boldsymbol{(0,0,e^{15} + e^{64},0,0,0)}$}

When $S$ corresponds to~\eqref{12+34}, we set
\begin{displaymath}
     \left\{
       \begin{array} {l}
         e^i = e^{\tfrac 34 m t} dx_i,\quad i = 1, 4, 5, 6,\\[3pt]
         e^2 = e^{m t} dx_2,\\[3pt]
         e^3 = -\tfrac 12 me^{\tfrac 32 mt}(\tfrac 32
dx_3+x_5dx_1+x_4dx_6),\\[3pt]
         e^7 = dt,
       \end{array}
     \right.
\end{displaymath}
whence
\begin{multline}
     \label{eq:metric12+34}
       \phantom{MMM} g = dx_2^2 +
       e^{-\tfrac 12 m t} ( dx_1^2 + dx_4^2 + dx_5^2 +
       dx_6^2)+\\[3pt]
        \tfrac 9{16} m^2 e^{m t} (dx_3 + \tfrac 23 x_5 dx_1 +
       \tfrac 23 x_4 dx_6)^2 + e^{-2m t} dt^2\phantom{MMMM}
\end{multline}
has holonomy $\Lie{SU}(3)\subset G_2$.  Restricting ourselves to $\lan x_2
\ran^\perp$, we obtain a metric on the product of a principal $T^1$-bundle
over $T^4$ with $\R$.

\medbreak
\subsection{The algebra $\boldsymbol{(0,e^{45},e^{64}+e^{51},0,0,0)}$}

Let us look at~\eqref{12,13+24} now:
\begin{displaymath}
     \left\{
       \begin{array} {l}
         e^i = e^{\tfrac 45 m t} dx_i,\quad i = 1, 6,\\
         e^2 = -\tfrac 35 m\,e^{\tfrac 65 m t}( dx_2 + \tfrac 23 x_4
dx_5),\\[3pt]
         e^3 = -\tfrac 35 m\,e^{\tfrac 75 m t}( dx_3 - \tfrac 23  x_1 dx_5
         + \tfrac 23  x_4 dx_6),\\[3pt]
         e^l = e^{\tfrac 35 m t} dx_l,\quad l = 4, 5,\\
         e^7 = dt
       \end{array}
     \right.
\end{displaymath}
allow to write down a previously unknown exceptional metric
\begin{proposition}
     \label{eq:new}
     Let $\n$ be the nilpotent Lie algebra defined by
     \begin{displaymath}
       e_2=[e_5,e_4], \q e_3=[e_6,e_4]=[e_1,e_5].
     \end{displaymath}
     Then with the above conventions, the solvmanifold $S$ relative to
     $\s=\n\oplus \R e_7$ carries a Riemannian metric
     \begin{multline}
       \label{eq:G2metric}
       g =   e^{-2m t} dt^2 + e^{-\tfrac 25 m t}(dx_1^2 + dx_6^2) +
       e^{-\tfrac 45 m t} (dx_4^2 + dx_5^2) + \\[3pt]
       \tfrac 9{25}m^2e ^{\tfrac 45 mt}(dx_3-\tfrac 23 x_1dx_5 + \tfrac 23 x_4
       dx_6)^2 + \tfrac 9{25}m^2e^{\tfrac 25 mt}(dx_2+\tfrac23 x_4dx_5)^2
     \end{multline}
     whose holonomy group is precisely $G_2$.\qed
\end{proposition}
\noindent
In sufficiently small neighbourhoods, this metric is clearly isometric to
\begin{multline*}
    \label{eq:isometricnewmetric}
    \phantom{MMM} ds^2= V^3dy^2 + V(dz_1^2 + dz_6^2) +
                        V^2(dz_4^2 + dz_5^2) + \\[3pt]
                        V^{-2}\bigl(dz_3+k(-z_1dz_5+z_4dz_6)\bigr)^2+
                        V^2(dz_2+k z_4dz_5)^2,\phantom{MMM}
\end{multline*}
with $V=ky$, on the product of $\R$ with a $T^2$-bundle over a $T^4$. The fibre
coordinates are $z_2,z_3$, whilst $y$ accounts for the $\R$ factor. The metric
\eqref{eq:G2metric} has a symmetry generated by the homothetic
Killing field
\begin{displaymath}
       Z = -\tfrac 5m\tfrac{\partial}{\partial
t}+4x_1\tfrac{\partial}{\partial x_1}+
        4x_6\tfrac{\partial}{\partial x_6}+3x_4\tfrac{\partial}{\partial x_4}+
        3x_5\tfrac{\partial}{\partial x_5}+
        \tfrac {21}5 mx_3\tfrac{\partial}{\partial x_3}+
        \tfrac {18}5 mx_2\tfrac{\partial}{\partial x_2},
\end{displaymath}
found by imposing invariance under a suitable scaling factor.
For appropriate Killing vector fields, this feature is common to all other
metrics in this section \cite{Gibbons-LPS:domain-walls}. Note that $Z$ does
not correspond to $e_7$, for
\begin{displaymath}
    (d Z^\flat)(\tfrac{\partial}{\partial t}, \tfrac{\partial}{\partial x_1}) =
    \tfrac{\partial}{\partial t} Z^\flat(\tfrac{\partial}{\partial x_1}) -
    \tfrac{\partial}{\partial x_1} Z^\flat(\tfrac{\partial}{\partial t}) -
    Z^\flat([\tfrac{\partial}{\partial t}, \tfrac{\partial}{\partial
x_1}])\neq 0.
   \end{displaymath}

\medbreak
\subsection{The algebra $\boldsymbol{(0,e^{45},e^{46},0,0,0)}$}

The coordinates
\begin{displaymath}
     \left\{
       \begin{array} {lcl}
         e^1 = e^{m t} dx_1, & & e^4 = e^{\tfrac 12 m t} dx_4,\\[3pt]
         e^2 = \tfrac 12 me^{\tfrac 54 m t}(-\tfrac 32 dx_2 + x_5 dx_4), & &
         e^i = e^{\tfrac 34 m t} dx_i,\quad i = 5, 6,\\[3pt]
         e^3 = \tfrac 12 me^{\tfrac 54 m t} (-\tfrac 32 dx_3 + x_6 dx_4), & &
         e^7 = dt,\\[3pt]
       \end{array}
     \right.
\end{displaymath}
relative to~\eqref{12,13} produce the metric
\begin{multline}
     \label{eq:metric12,13}
       \phantom{M} g = dx_1^2 + \tfrac 9{16}m^2
       e^ {\tfrac 12 mt} (dx_2 - \tfrac 23 x_5 dx_4)^2
       + \tfrac 9{16}m^2 e^{\tfrac 12 m t} (dx_3 - \tfrac 23 x_6 dx_4)^2\\[3pt]
       + e^{-m t} dx_4^2 + e^{-\tfrac 12 m t} (dx_5^2 + dx_6^2) + 
e^{-2m t} dt^2.
\end{multline}
It has holonomy $\SU(3)$, too. This can be recovered by looking at $\lan
x_1\ran^\perp=\R^6$, and the induced metric on the product of $\R$ with a
principal $T^2$-bundle over $T^3$.

\medbreak
\subsection{The Iwasawa algebra}

\noindent
The algebra~\eqref{13+42,12+34} comes equipped with 1-forms
\begin{equation}\label{iwabasis}
     \left\{
       \begin{array} {l}
         e^i = e^{\tfrac 23 m t} d x_i,\quad i = 1, 4, 5, 6,\\[3pt]
         e^2 =\tfrac m3 e^{\tfrac 43 m t} (2dx_2 +x_6 dx_1 - x_4 dx_5),\\[3pt]
         e^3 =-\tfrac m3 e^{\tfrac 43 m t} (2dx_3 + x_5 dx_1 + x_4 dx_6),\\
         e^7 = dt.
       \end{array}
     \right.
\end{equation}
Then
\begin{equation}
\label{eq:iwametric}
\begin{array}{r}
       g =
      \tfrac 49 m^2 e^{\tfrac 23 m t} (d x_3 +
       \tfrac
12 x_5 dx_1 + \tfrac 12 x_4 dx_6)^2 +
       \tfrac 49 m^2 e^{\tfrac
23 m t}  (d x_2 +\tfrac 12 x_6 dx_1 -
      \tfrac 12 x_4 dx_5)^2 +
\\[3pt]
      e^{-\tfrac 23 mt} (dx_1^2 + dx_4^2 + dx_5^2 + dx_6^2) +
e^{-2mt}dt^2
      \phantom{MMMMMM}
\end{array}
\end{equation}
has
holonomy $G_2$.

\medbreak
\subsection{The algebra $\boldsymbol{(0,0,e^{15},e^{25},0,e^{12})}$}

Eventually,
\begin{multline}
\label{eq:instantonmetric}
   g =
   e^{-\tfrac 45 mt} (dx_1^2 + dx_2^2
+ dx_5^2) +
   \tfrac 9{25}m^2e^{\tfrac 25 m t}(d x_3 + \tfrac 23 x_5
dx_1)^2\\[3pt]
   + \tfrac 9{25} m^2 e^{\tfrac 25 mt}  (d x_4 - \tfrac
23 x_2
   dx_5)^2 + \tfrac 9{25}m^2 e^{\tfrac 25 m t} (d
   x_6 +
\tfrac 23 x_2 dx_1)^2 + e^{-2m t} dt^2
\end{multline}
belongs
to~\eqref{12,13,23}, by means of
\begin{displaymath}
     \left\{
       \begin{array} {l}
         e^i = e^{\tfrac 35 m t} d x_i,\qquad i = 1, 2, 5,\\[3pt]
         e^3 = -\tfrac 15 me^{\tfrac 65 m t}(3dx_3 +2x_5 dx_1),\\[3pt]
         e^4 = -\tfrac 15 me^{\tfrac 65 m t}(3dx_4 -2x_2 dx_5),\\[3pt]
         e^6 = -\tfrac 15 me^{\tfrac  65 m t}(3dx_6 +2x_2 dx_1),\\
         e^7 = dt.
       \end{array}
     \right.
\end{displaymath}
These expressions identify $g$ as a
$G_2$ holonomy metric on $\R$ times a $T^3$-bundle over $T^3$.

\bigbreak\noindent
Leaving the torus' flat structure aside, three of the metrics found
have reduced
holonomy, i.e.~$\SU(2)$ and $\SU(3)$. These are attached to the algebras of
Theorem~\ref{theorem} containing an Abelian summand
\begin{displaymath}
    \R^3\oplus \h_3,\qq \R \oplus \h',\qq \R\oplus \h''
\end{displaymath}
for given $\h',\h''$. This hints that the $G_2$ metrics
could be reduced to lower dimensional structures of special type, with the
same philosophy pursued in \cite{Apostolov-S:K-G2}. On the other hand
\eqref{eq:G2metric}, \eqref{eq:iwametric} and \eqref{eq:instantonmetric}
are proper holonomy $G_2$ metrics, and are indeed built from algebras with
an irreducible and more complicated structure.  All metrics are
scale-invariant.  This is because the corresponding groups $S$ can be
decomposed into irreducible de~Rham factors which are scale-invariant.  The
correspondence between $\n$ and the holonomy of the Ricci-flat metric $g$
supported by its rank-one solvable extension $\s$ is summarised in the
table.
\begin{table}[!ht]
   \begin{displaymath}
     \begin{array}{|c|c|}
      \hline \vs \hbox{Nilpotent algebra $\n$} & \hbox{Holonomy} \\
      \hline\hline \vs (0,0,e^{15},0,0,0) & \SU(2) \\
      \hline \vs       (0,0,e^{15} + e^{64},0,0,0) & \SU(3)\\
      \hline \vs       (0,e^{45},e^{64}+e^{51},0,0,0) & G_2\\
      \hline \vs       (0,e^{45},e^{46},0,0,0) & \SU(3)\\
      \hline \vs       (0,e^{16}+e^{45},e^{15}+e^{64},0,0,0) & G_2\\
      \hline \vs       (0,0,e^{15},e^{25},0,e^{12}) & G_2\\
      \hline
     \end{array}
   \end{displaymath}
\end{table}

\section{Evolving the nilpotent $\SU(3)$ structure}\label{sec:evolution}

\noindent
Let us turn to $\Gamma \backslash N \times \R$ and consider on the
nilmanifold $\Gamma \backslash N$ the $\SU(3)$ structure induced by that of
$N$. In this section we wish to explain how one can use the so-called
evolution equations discovered in~\cite{Hitchin:forms}. These predict the
deformation in time of special kinds of $\SU(3)$ structures and their
ability to give rise to metrics with holonomy contained in $G_2$.  As a
matter of fact, it turns out that half-flat $\SU(3)$-manifolds represent
the natural class with the potential to evolve along the flow of the
differential system and be preserved by it at the same time.  We thus
assume that $\omega(\t),\psp(\t)$ is an $SU(3)$ structure depending on a
locally defined real parameter $\t\in \R$. We may then regard the resulting
$7$-manifold as fibring over an interval, which accounts for a `dynamic'
inclusion of $\SU(3)$ in the exceptional group.  The central point is that
the fundamental forms evolve according to the differential equations
\begin{equation}
    \label{eq:evolution}
    \left\{
\begin{array}{l}
        \displaystyle \hd\omega  =  \p\psp,\\
  \displaystyle \hd\psm  =  -\omega\w\p\omega.
      \end{array}
  \right.
\end{equation}
The compatibility relations restraining the
almost Hermitian structure
\begin{gather}
\omega\w\psp=0,\label{eq:compatibility1}\\
    \psp\wedge\psm=\tfrac 23
\omega^3\label{eq:compatibility2}
\end{gather}
are preserved in time.  Considering the general difficulty in solving
system \eqref{eq:evolution}, half-flatness provides the simplest examples
of evolution structures --- other than nearly K\"ahler ones. Besides, as
all data in question is analytic, the solution is uniquely determined,
hence one can expect the outcome to resemble one of the metrics of the
previous section
\begin{proposition}
\label{6=5}
  Any of the Ricci-flat metrics on the solvable Lie group 
  $S$ with structure equations \eqref{h3} -- \eqref{12,13,23}, can be
  obtained evolving the $SU(3)$ structure on the $2$-step nilmanifold
  $\Gamma \backslash N$.
\end{proposition}
\noindent
For understandable reasons the discussion will omit the case of $T^6$, for
which the results of this section hold anyway, if trivially.  To avoid
being too long, the proof will rely on the detailed description of the
technique for $\h_3\oplus \R^3$ and the Iwasawa algebra only.  The forms
$\psp(0),\ \oo(0)$ will be indicated by $\psp_0,\ \oo_0$.

\medbreak\noindent {\bf First example.}  We begin by considering the
nilpotent Lie algebra
\begin{displaymath}
(0,0,\tfrac 23 m e^{15},0,0,0)
\end{displaymath}
underlying that of~\eqref{h3}.  The forms defining the structure on $\Gamma
\backslash N$ are deformed by means of exact elements in
Chevalley--Eilenberg's cohomology
\begin{displaymath}
   \Omega^3_{\textrm{exact}}=\lan 
e^{125},e^{145},e^{156}\ran, \quad
   \Omega^4_{\textrm{exact}}=\lan 
e^{1245},e^{1256},e^{1456}\ran.
\end{displaymath}
We introduce the deformation functions, all depending upon $\t$
\begin{displaymath}
\left\{
\begin{array}{l}
      \tfrac 12 \oo(\t)=P(\t)\tfrac 12\oo_0+
D(\t)e^{1245}+E(\t)e^{1256}+F(\t)e^{1456},\\[4pt]
\psp(\t) = 
Q(\t)\psp_0+A(\t)e^{125}+B(\t)e^{145}+C(\t)e^{156}.
\end{array}
\right.
\end{displaymath}
By asking $P(0)=Q(0)=1,\ D(0)=E(0)=F(0)=A(0)=B(0)=C(0)=0$, one is able to
regain the initial structure~\eqref{eq:standardstructure} at time $\t=0$.
The expression for $\oo(\t)$ suggests that the K\"ahler form, uniquely
determined up to sign, must be of the following kind
\begin{displaymath}
   \omega 
(\t)=x(\t)e^{14}+y(\t)e^{23}+z(\t)e^{56}+w(\t)e^{25}+j(\t)e^{12},
\end{displaymath}
for certain functions satisfying $x(0)=z(0)=-y(0)=1,\ w(0)=j(0)=0$.  In the
following, the explicit dependence upon $\t$ will be dropped, with the
convention that the relations hold for all appropriate values of time. By
computing $\omega\wedge\omega$ one finds the relations
\begin{equation}
   \label{rel}
   zx=P+F,\quad 
xy=yz=-P,\quad xw=-D,\quad zj=E.
\end{equation}
In addition, the primitivity of $\pspm$ underlying
equation~\eqref{eq:compatibility1} implies $yB=jQ,\ yC=wQ$.  The first
differential equation of~\eqref{eq:evolution} compares
\begin{displaymath}
  \textstyle\p 
\psp=(Q'+A')e^{125} +Q'(-e^{345}+e^{136}+e^{246})+ 
B'e^{145}+C'e^{156}
\end{displaymath}
with 
\begin{displaymath}
\hd\omega=\tfrac 23 my 
e^{125},
\end{displaymath}
dashed letters denoting derivatives with respect to $\t$. This gives
$Q(\t)=1, B(\t)=C(\t)=0$ for all $\t$'s, and $A'(\t)=\tfrac 23 my(\t)$.
Therefore
\begin{displaymath}
\psp=(A+1)e^{125}-e^{345}+e^{136}+e^{246}. 
\end{displaymath}
From~\eqref{rel} one has $j(\t)=w(\t)=0$, hence $D(\t)=E(\t)=0$.  The
determination of an orthonormal basis of 1-forms enables one to get
$\psm(\t)$. Inspired by equations~\eqref{lambda coefficients}, one defines
\begin{displaymath}
  \lambda^a\1,\ \lambda^b\2,\ \lambda^c\3,
  \ 
\lambda^b\4,\ \lambda^a\5,\ \lambda^b\6,
\end{displaymath}
for some non-zero function $\lambda=\lambda(\t)$ with $\lambda(0)=1$.  For
$\psp=\lambda^{2a+b}
e^{125}+\lambda^{a+b+c}(-e^{345}+e^{136})+\lambda^{3b}e^{246}$ to resemble
the previous expression one takes $a+b+c=0=3b$. One of many possible
choices is $a=\tfrac 12=-c$, so that now $\pspm$ assume the form
\begin{displaymath}
\psp=\lambda 
e^{125}-e^{345}+e^{136}+e^{246},\quad 
  \psm=\lambda^{1/2}(e^{126}- 
e^{245}-e^{135})-\lambda^{-1/2}e^{346}
\end{displaymath}
and thus $\lambda(\t)=A(\t)+1$. The second evolution equation
\begin{displaymath}
  -(P'+F')e^{1456} 
+P'(e^{1423}+e^{2356}) 
= -\omega\omega' = 
  \hd\psm=\tfrac 
{2m}{3\sqrt\lambda} 
e^{1456}
\end{displaymath}
implies $P(\t)=1$ and $F'=-\tfrac{2m}{3\sqrt{A+1}}$.  Eventually, volume
normalisation \eqref{eq:compatibility2} says that
\begin{displaymath}
\left\{
\begin{array}{l}
      \sqrt{A+1}=-\tfrac 23 m (A')^{-1} \\
      A(0)=0.
    \end{array}
  \right.  
\end{displaymath}
The solution to this initial value problem reads
\begin{displaymath}
A(t)= 
(1-m\t)^{2/3}-1,
\end{displaymath}
so the geometric structure is evolving according to
\begin{gather*}
\psp(\t)=\sqrt[3]{(1-m\t)^2}\,e^{125}-e^{345}+e^{136}+e^{246},\\
\omega(\t)=\sqrt[3]{1-m\t}\,(e^{14}+e^{56})-\frac{1}{\sqrt[3]{1-m\t}}\,e^{23}.
\end{gather*}
The associated metric
\begin{displaymath}
  g = 
(1-m\t)^{2/3}\bigl((\1)^2+(\5)^2\bigr)+\sum_{i=2,4,6}(e^i)^2+
(1-m\t)^{-2/3}(\3)^2+d\t^2
\end{displaymath}
mirrors precisely \eqref{7metric}: the appropriate coordinate system
$\{x_i\}$ on $\R^6$ is given by
\begin{displaymath}
  e^i=dx_i,\q i\neq 3, \qq 
\3=-\tfrac 23 m(dx_3+x_5dx_1),
\end{displaymath}
and the correspondence follows once one identifies $1-m\t$ with the
conformal factor $e^{-mt}$.

\bigbreak\noindent {\bf Second example.}  Let us look at the Iwasawa Lie
algebra
\begin{displaymath}
\bigl(0,-\tfrac 
13 
m(e^{16}+e^{45}),\tfrac 13 
m(e^{15}-e^{46}),0,0,0\bigr),
\end{displaymath}
whose evolution also deserves a detailed description. To start with, the
deformation is given by
\begin{displaymath}
\begin{array}{l}
\psp=Q\psp_0+Ae^{145}+Be^{146}+Ce^{156}+D^{456},\\[4pt]
\tfrac 12 
\oo=\tfrac 12 
  P\oo_0 + Fe^{1456} + G(e^{1435}-e^{1426}) 
+ 
H(-e^{1436}-e^{1425}) +\\[4pt]
  \phantom{MMMMMM} 
L(e^{1356}+e^{2456}) + 
M(e^{3456}-e^{1256}),
\end{array}
\end{displaymath}
the latter telling that the K\"ahler form is
\begin{displaymath}
\omega=xe^{14}+ye^{23}+ze^{56}+u(e^{35}-e^{26})+v(-e^{36}-e^{25})+
\rho(e^{13}+e^{24})+\sigma(e^{34}-e^{12}).
\end{displaymath}
This first relations obtained by comparison are
\begin{equation}
\label{eq:xyz}
  \begin{array}{c}
    -P=xy=yz, 
\quad P+F=xz,\\
    G=xu, \quad H=xv, \quad L=z\rho, \quad 
M=z\sigma,\\
    u\rho+v\sigma=0, \quad 
-u\sigma+v\rho=0.
\end{array}
\end{equation}
Then the evolution of $\psp$ immediately annihilates $A, B, C, D$ and
yields
\begin{displaymath}
  Q'=\tfrac 13 my,
\end{displaymath}
as $\psp_0$ is exact. On the other hand $0=\tfrac 1Q
\omega\psp=-2e^{23}(ue^{456}+ve^{156}+\rho e^{146}-\sigma e^{145})$ forces
the vanishing of most of the remaining coefficients
\begin{displaymath}
  \psp=Q\psp_0, \qquad \tfrac 12 \oo=\tfrac 
12 P\oo_0 + Fe^{1456}.
\end{displaymath}
The orthonormal basis $\lambda^a e^i,\ 
i=1,4,5,6,\ \lambda^b e^j,\ j=2,3\,$ is chosen in order to mimic equations
\eqref{iwabasis}. The actual values of $a, b\in\R$ are irrelevant at
present, one could for example take $a=-1,b=1$.  In any case
\begin{displaymath}
\psp=\lambda^{2a+b}\psp_0,\q 
\psm=\lambda^{2a+b}\psm_0, \ 
\textrm{with}\ 
  \lambda^{2a+b}=Q. 
\end{displaymath}
From $\hd\psm=-\omega\omega'$ one obtains $P=1$ and $F'=\tfrac 43 mQ$, so
that each of $\oo,\psp$ evolves in one direction only. This could have been
predicted by counting dimensions, see~\cite{Ketsetzis-S:Iwasawa}. The first
line in~\eqref{eq:xyz} tells that $x=z=-1/y=\sqrt{F+1}$, and
\eqref{eq:compatibility2} gives the quartic curve
\begin{displaymath}
Q^2=\sqrt{F+1},
\end{displaymath}
so the evolution equations are equivalent to the first order system
\begin{displaymath}
\left\{\displaystyle
     \begin{array}{l}
     Q'(\t)= 
-\frac{m}{3}Q^{-2}\\[6pt]
     F'(\t)=-\tfrac 43 m(F+1)^{1/4}\\[5pt]
Q(0)=1, \q F(0)=0.
     \end{array}
   \right.
\end{displaymath}
The first equation is solved by $Q(\t)=(1-m\t)^{1/3}$, hence
$F(\t)=(1-m\t)^{4/3}-1$. Eventually, the $\SU(3)$ structure results in
\begin{gather*}
\psp(\t)=(1-m\t)^{1/3}\,(e^{125}-e^{345}+e^{136}+e^{246}),\\
\omega(\t)=(1-m\t)^{2/3}\,(e^{14}+e^{56})-\dfrac{1}{(1-m\t)^{2/3}}\,e^{23}.
\end{gather*}
These data identify the non-integrable complex structure $-J_3$ studied in
a broader context by \cite{Abbena-GS:aH-nil} and the flow corresponds to
the one given in \cite[ex.~2 (iii)]{Apostolov-S:K-G2}. A glance at the
level curves of the background Hamiltonian function confirms that the
almost complex structure degenerates at time $\t=1/m$.  The corresponding
$G_2$-metric
\begin{displaymath}
\displaystyle
g=(1-m\t)^{2/3}\,\sum_{\substack{i=1,4,\\5,6}}(e^i)^2 +
\frac{1}{(1-m\t)^{2/3}}\,\left( (\2)^2 +(\3)^2\right) + 
d\t^2
\end{displaymath}
has its counterpart in \eqref{eq:iwametric} when $t=\tfrac 1m \ln |1-m\t|$.

\bigbreak\noindent {\bf The other cases.}  ${\bf 1.}$ With the same
technique we tackle the 6-dimensional Lie algebra with non-trivial brackets
$[e_5,e_1]=\tfrac 12 m e_3=[e_4,e_6]$, isomorphic to
$(0,0,0,0,0,e^{12}+e^{34})$.  As time goes by, the outcoming $\SU(3)$
geometry is described by
\begin{gather*}
\psp(\t)=\sqrt{1-m\t}\,(e^{125}+e^{246})+e^{136}-e^{345},\\
\omega(\t)=\sqrt{1-m\t}\,(e^{14}+e^{56})-\frac{1}{\sqrt{1-m\t}}e^{23},
\end{gather*}
so the Riemannian structure is
\begin{displaymath}
g=\sqrt{1-m\t}\sum_{\substack{j=1,4,\\5,6}}(e^j)^2+(\2)^2+
\frac{1}{1-m\t}(\3)^2+d\t^2,
\end{displaymath}
in agreement with \eqref{eq:metric12+34}.

\medbreak ${\bf 2.}$ We next apply the evolution machinery to
$\bigl(0,-\tfrac 25 me^{45}, \tfrac 25 m(e^{15}-e^{46}),0,0,0\bigr)$. The
four-form flows according to $\tfrac 12 \oo(\t)=
(1-m\t)^{6/5}e^{1456}-e^{1423}-e^{2356}$, whose square root provides
\begin{displaymath}
\omega(\t)=(1-m\t)^{3/5}(e^{14}+e^{56})-(1-m\t)^{-3/5}e^{23}, 
\end{displaymath}
\noindent while the $3$-form is 
\begin{displaymath}
     \psp(\t) = 
(1-m\t)^{2/5}(-e^{345}+e^{125}+e^{246})+e^{136}.
\end{displaymath}
The Riemannian metric
\begin{displaymath}
   \begin{split}
     g & = 
(1-m\t)^{2/5}\,\bigl((\1)^2+(\6)^2\bigr)+
(1-m\t)^{4/5}\,\bigl((\4)^2+(\5)^2\bigr)\\
& \quad 
+(1-m\t)^{-2/5}\,(\2)^2+(1-m\t)^{-4/5}\,(\3)^2+d\t^2,
\end{split}
\end{displaymath}
basically recovers that of Proposition~\ref{eq:new}.

\medbreak ${\bf 3.}$ The deformation of the Lie structure corresponding to
$d\2=-\tfrac 15 me^{45},\, d\3=-\tfrac 15 me^{46}$ is
\begin{gather*}
\psp=e^{125}+e^{136}+\sqrt{1-m\t}\,(e^{246}-e^{345}),\\
\omega=\sqrt{1-m\t}\,(e^{14}+e^{56})-\frac{1}{\sqrt{1-m\t}}e^{23},
\end{gather*}
hence we get
\begin{displaymath}
g=(\1)^2+\frac{1}{\sqrt{1-m\t}}\bigl((\2)^2+(\3)^2\bigr)+
(1-m\t)\,(\4)^2+\sqrt{1-m\t}\,\bigl((\5)^2+(\6)^2\bigr)+d\t^2,
\end{displaymath}
related to \eqref{eq:metric12,13}.

\medbreak ${\bf 4.}$ At last, let us write what happens to $(0,0,\tfrac 25
me^{15}, \tfrac 25 me^{25},0,\tfrac 25 me^{12})$. After some computations
one finds that
\begin{gather*}
  \psp=(1- m\t)^{6/5}\,e^{125}+e^{136}+e^{246}-e^{345},\\
  \omega=(1-m\t)^{1/5}\,(e^{14}-e^{23}+e^{56})
  \end{gather*}
  are compatible with the metric
\begin{displaymath}
g=(1-m\t)^{4/5}\sum_{i=1,2,5}(e^i)^2+ 
(1-m\t)^{-2/5}\sum_{j=3,4,6}(e^j)^2+d\t^2,
\end{displaymath}
see 
\eqref{eq:instantonmetric}.\qed

\medbreak It is no coincidence that in all cases the
identification between the `evolved' holonomy metrics and the ones found in
\S\ref{sec:metrics} is attained by uniformly putting $\t = 1/m(1 - e^{-
  mt})$ and using global coordinates $x_1,\ldots, x_6$ on the nilpotent Lie
group $N$ to represent the left-invariant forms $\{e^i\}$. Which brings to
the completeness' properties of the metrics. We are always in presence of a
unique singularity, determined by $f(t)=\exp(-mt)$ or, if one prefers, by
the linear function $\tilde f(\t)= 1-m\t$. This means that away from the
degeneration, all metrics are complete in one direction of time, i.e.~the
tensors $g, \psp, \omega$ describe a smooth
structure for $\t\in(-\infty,\t_0],\ \t_0<1/m$, see related discussion in
\cite{Apostolov-S:K-G2}.

\bibliographystyle{amsplain}

\providecommand{\bysame}{\leavevmode\hbox 
to3em{\hrulefill}\thinspace}
\providecommand{\MR}{\relax\ifhmode\unskip\space\fi MR }
\providecommand{\MRhref}[2]{%
\href{http://www.ams.org/mathscinet-getitem?mr=#1}{#2}
}
\providecommand{\href}[2]{#2}

\end{document}